\documentclass[reqno,10pt]{amsart}
\usepackage{xspace}
\usepackage[bookmarksnumbered,colorlinks]{hyperref}
\usepackage{graphics}
\usepackage{fourier}
\usepackage[english]{babel}
\newcommand{\bt}{\begin{theorem}}                     
\newcommand{\et}{\end{theorem}}                       
\newcommand{\bd}{\begin{definition}}                  
\newcommand{\ed}{\end{definition}}                    
\newcommand{\bl}{\begin{lemma}}                       
\newcommand{\el}{\end{lemma}}                                   
\newcommand{\bpr}{\begin{proposition}}                  
\newcommand{\epr}{\end{proposition}}                    
\newcommand{\bere}{\begin{remark}}                      
\newcommand{\ere}{\end{remark}}                         
\newcommand{\beq}{\begin{equation}}
\newcommand{\eeq}{\end{equation}}
\def\bal#1\eal{\begin{align}#1\end{align}}              
\def\baln#1\ealn{\begin{align*}#1\end{align*}}          
\def\bml#1\eml{\begin{multline}#1\end{multline}}        
\def\bmln#1\emln{\begin{multline*}#1\end{multline*}}  
\def\bga#1\ega{\begin{gather}#1\end{gather}}
\def\bgan#1\egan{\begin{gather*}#1\end{gather*}}
\newcommand{\de}{\mathrm{d}}                        
                     
\newcommand{\N}{\ensuremath{\mathbb{N}}\xspace}     
\newcommand{\R}{\ensuremath{\mathbb{R}}\xspace}     
\newcommand{\eps}{\varepsilon}                      
\newcommand{\To}{\longrightarrow}                   
\newcommand{\inte}{\int_0^1\!\!}

\newtheorem{theorem}{Theorem}[section]

\newtheorem{lemma}[theorem]{Lemma}
\newtheorem{proposition}[theorem]{Proposition}
\theoremstyle{definition}
\newtheorem{definition}[theorem]{Definition}
\theoremstyle{remark}
\newtheorem{remark}[theorem]{Remark}
\newcommand{\dist}{\ensuremath{\mathrm{dist}}\xspace}
\newcommand{\nablag}{\ensuremath{\nabla^h}\xspace}

\title[Morse theory of causal geodesics in a stationary spacetime]%
{Morse theory of causal geodesics in a stationary spacetime via Morse theory of geodesics of a Finsler metric}

\author[E. Caponio]{Erasmo Caponio}
\address{Dipartimento di Matematica, Politecnico di Bari, Via Orabona 4,
70125, Bari, Italy}
\email{caponio@poliba.it}
\thanks{EC and AM are supported by M.I.U.R. Research project PRIN07 ``Metodi Variazionali e 
Topologici nello Studio di Fenomeni Nonlineari''
}

\author[M. A. Javaloyes]{Miguel \'Angel Javaloyes}
\address{Departamento de Geometr\'{\i}a y Topolog\'{\i}a.
  Facultad de Ciencias, Universidad de Granada.
 Campus Fuentenueva s/n, 18071 Granada, Spain}
\email{majava@ugr.es, ma.javaloyes@gmail.com}
\thanks{MAJ is partially supported by Regional J.
Andaluc\'{\i}a Grant P06-FQM-01951, by Fundaci\'on S\'eneca project
04540/GERM/06   and by Spanish MEC Grant MTM2007-64504}

\author[A. Masiello]{Antonio Masiello}
\address{Dipartimento di Matematica,
Politecnico di Bari, Via Orabona 4,
70125, Bari, Italy}
\email{masiello@poliba.it}

\subjclass[2000]{53C22, 53C50, 53C60, 58E05}

\keywords{stationary Lorentzian manifolds, light rays, Morse theory, conjugate points,  Finsler metrics.}

\date{}

\begin{document}

\begin{abstract}
We show  that the index of  a lightlike  geodesic in a conformally standard stationary spacetime $(\mathcal M_0\times\R,g)$ is equal to the index of its spatial projection as a geodesic of a Finsler metric $F$ on $\mathcal M_0$ associated to $(\mathcal M_0\times\R,g)$. Moreover we obtain the Morse relations of lightlike geodesics connecting a point $p$ to a curve  $\gamma(s)=(q_0,s)$  by using Morse theory on the Finsler manifold $(\mathcal M_0,F)$. To this end, we prove a splitting lemma for the energy functional of a Finsler metric.  Finally, we show  that the reduction to Morse theory of a Finsler manifold  can be done  also for timelike geodesics.

\end{abstract}
\maketitle
\begin{section}{Introduction}\label{intro}
Since the seminal paper  \cite{Uhlenb75},  Morse theory has been applied successfully to  spacetime geometry (Lorentzian manifolds) and global problems of  general relativity. For instance, a consequence of the Morse relations of lightlike geodesics proved in \cite{Uhlenb75} is  that  on a  contractible globally hyperbolic spacetime,
whose metric  satisfies a suitable  growth condition, the number of images produced by a gravitational lens
is odd (or infinity) \cite{McKenz85}. Gravitational lensing is the phenomenon where  the gravitational field of a  galaxy, located between an observer  and a star,  bends the light rays  emitted by the star and   focuses them at the same instant of observation, causing the observer to see multiple  images of the same star (see e. g. \cite{ScEhFa92}).

After the papers \cite{Uhlenb75,McKenz85}, Morse theory has been applied to compute the number of lightlike geodesics between an event and a timelike curve, on different classes of spacetimes and for different types of lenses and sources (see e. g. \cite{FoGiMa95,GiaLom99,GiMaPi02,HasPer06,Petter92,Petter95}).

We recall that a Lorentzian manifold $(\mathcal M,g)$ is a smooth  connected manifold $\mathcal M$ endowed with a symmetric non-degenerate tensor field $g$ of type $(0,2)$ having index $1$.
The geodesics of $(\mathcal M,g)$ are the critical points of the energy functional of the metric $g$
\beq
z\mapsto\tfrac{1}{2}\inte g(z)[\dot z,\dot z]\de s.\label{nrglorentzian}\eeq
So they are the smooth curves $z\colon[a,b]\to\mathcal M$ satisfying the equation $\nabla_{\dot z}\dot z=0$, where $\nabla$ is the  Levi-Civita connection of the metric $g$.
If $z$ is a geodesic,  the function $s\mapsto g(z(s))[\dot z(s),\dot z(s)]\colon=E_{z}$  is constant. According to the sign of $E_{z}$, a  geodesic is said {\em timelike} if
$E_{\gamma}<0$, {\em lightlike} if
$E_{\gamma}=0$ and {\em spacelike} if
$E_{\gamma}>0$ or $\dot z=0$. Timelike and lightlike geodesics are also called {\em causal}.
Such a terminology is also used for any vector in any tangent space and for any piecewise smooth curve iff its  tangent vector field has the same character at any point where it is defined.

A striking difference with the Riemannian case is  that the energy functional of a Lorentzian metric is unbounded both from below and above and the Morse index of its critical points is
$+\infty$. The common strategy used to develop Morse theory of a Lorentzian manifold is to consider only a particular type of  geodesics (timelike or lightlike), to restrict the index form
to the vector fields that are orthogonal to the geodesic, in the timelike case, and  orthogonal modulo the vector fields pointwise collinear to the velocity vector field of the geodesic,
in the lightlike case, and to use the length functional on a finite dimensional approximation of the  path space (see \cite[Ch. 10]{BeErEa96} and the references therein).
Another approach is to substitute  the energy functional with a functional which has nice variational properties.  This works for lightlike geodesics, which are the critical points of the {\em arrival time functional} (see \cite{Uhlenb75}),  or for particular kinds of Lorentzian manifolds as the  standard stationary ones (see \cite{BenMas92,Masiel94}).\footnote{For recent results about the Morse index theorem in the spacelike case and  the Morse relations  for all type of geodesics see respectively \cite{PiTa02} and \cite{AbbMaj06}.}

The  aim of this paper is twofold: to show that for a  standard stationary Lorentzian manifold,  Morse theory of causal geodesics can be reduced to Morse theory of geodesics of a Finsler manifold of Randers type associated to the spacetime; to show that Morse theory for geodesics  connecting two points on a Finsler  manifold can be casted in a purely infinite dimensional setting without using finite dimensional approximations.

In regard to the first aim, we  show that the number of conjugate instants (counted with their multiplicity) along a lightlike [resp.  timelike] geodesic is equal to that of the corresponding Finslerian geodesic (Theorem~\ref{index}) [resp. Theorem~\ref{timelike}]. Moreover the Morse relations of lightlike [resp. timelike  parametrized with respect to the proper time on a
given interval] geodesics  joining a point with a timelike curve on the spacetime
can be obtained from the Morse relations of  the geodesics joining two points on the Finsler manifold (Theorem~\ref{Morseluce}) [resp.  Theorem~\ref{timelike}]. Although this reduction  is very natural and convenient, stationary spacetimes seem to be the only type of spacetimes where it works fine, without  leaving  the  realm of strongly convex Finsler metrics (cf. also \cite{GHWW08}).

We recall that a \emph{Finsler metric} $F$ on a manifold $M$ is a
continuous function $F\colon TM\to[0,+\infty)$  such that
\begin{itemize}
 \item $F$ is smooth  on $TM\setminus 0$;
\item $F$ is fiberwise positively  homogeneous of degree one, that is $F(x,\lambda y)=\lambda F(x,y)$, for all $x\in M$,   $y\in T_x M$ and  $\lambda>0$;
\item $F$ has fiberwise strongly convex square, that is
\[
{\mathbf{g}}_{ij}(x,y)=\left[\frac{1}{2}\frac{\partial^2 (F^2)}{\partial y^i\partial y^j}(x,y)\right]
\]
is positively defined for any $(x,y)\in TM\setminus 0$.
\end{itemize}
By the Euler's theorem we have that $F^2(x,y)=\mathbf{g}(x,y)[y,y]$.
A Finsler metric is said of Randers type if
\[F(x,y)=\sqrt{\alpha(x)[y,y]}+\omega(x)[y],\]
where $\alpha$ is a Riemannian metric on $M$ and $\omega$ is a $1$-form on $M$ having norm with respect to $\alpha$ strictly less than $1$ (see \cite[p. 17]{BaChSh00}).

The length of a piecewise smooth  curve
$\gamma\colon [a,b]\subset\R \to M$ with respect to the Finsler metric $F$  is defined
by $L(\gamma)=\int_a^b\!\! F(\gamma(s),\dot\gamma(s))\de s$. Thus the distance  between two arbitrary points $p,\ q\in M$ is given by
\beq\label{finsleriandist}
\dist(p,q)= \inf_{\gamma\in C(p,q)}L(\gamma),
\eeq
where $C(p,q)$ is the set of all piecewise smooth curves $\gamma\colon[a,b]\to M$ with $\gamma(a)=p$ and $\gamma(b)=q$.
The distance function \eqref{finsleriandist} is nonnegative and satisfies the triangle inequality, but it is not symmetric as $F$ is non-reversible. Thus one has  to distinguish the order  of   a pair  of points in $M$ when speaking about  distance. As a consequence, one is naturally led to the notions of forward and backward  Cauchy sequences and completeness (see \cite[\S 6.2]{BaChSh00}):
a sequence $\{x_n\}\subset M$ is called {\em forward} [resp. {\em backward}] {\em Cauchy sequence} if for all $\eps>0$ there exists $\nu\in\N$ such that, for all $\nu\leq i\leq j$, $\dist(x_i,x_j)\leq \eps$ [resp. $\dist(x_j,x_i)\leq \eps$]; $(M,F)$ is {\em forward complete} [resp. {\em backward complete}] if all forward [resp. backward] Cauchy  sequences converge.

The  geodesics $x:[0,1]\to M$  of a Finsler manifold $(M,F)$  parametrized with constant speed $F(x,\dot x)$
are the curves $x$ satisfying the equation
\[D_{\dot x}\dot x=0,\]
where $D_{\dot x}\dot x$ is the {\em Chern covariant derivative} of $\dot x$ along $x$ with reference vector $\dot x$ (see \cite[Chapter 5 and Exercise 5.2.5]{BaChSh00}). As it is shown for example in \cite[Proposition 2.3]{CaJaMa07a}, the geodesics parametrized with constant speed joining two given points $p_0,q_0\in {M}$ coincide with  the critical points of the energy functional
\[
E(x)=\frac{1}{2}\inte F^2(x,\dot x)\de s
\]
defined on the manifold $\Omega_{p_0,q_0}(M)$, which is
the collection  of the curves $x\colon[0,1]\rightarrow M$ such that $x(0)=p_0,\ x(1)=q_0$  and having $H^1$-regularity, that is $x$ is absolutely continuous and the integral $\inte h(x)[\dot x,\dot x]\de s$ is finite. Here $h$ is any complete Riemannian  metric on $M$.  It is well known that $\Omega_{p_0,q_0}(M)$ is a Hilbert manifold
modeled on any of the equivalent  Hilbert  spaces of $H^1$-sections, with vanishing endpoints, of the pulled back bundle $x^*TM$, $x$  any  regular curve in $M$ connecting $p_0$ to $q_0$
\cite[Proposition 2.4.1]{Klinge82}. The Riemannian metric on $\Omega_{p_0,q_0}(M)$ is given by
\[
\langle X,Y\rangle=\inte h(x)[\nabla^h_{\dot x} X,\nabla^h_{\dot
x}Y]\de s,
\]
for every  $H^1_0$-section, $X$ and $Y$ of $x^*T M$,
$\nablag$ being the Levi-Civita connection of the metric $h$.

As in Riemannian geometry, Jacobi vector fields are  the vector fields along the geodesic $x$ which give rise to variations of $x$ by means of  geodesics (parametrized with constant Finslerian speed), see \cite[\S 5.4]{BaChSh00} or \cite[\S 11.2]{Shen01}. A conjugate instant $\bar s$ along the geodesic $x\colon[0,1]\to M$ is a value of the parameter $s$ such that there exists a Jacobi vector field $J$, with $J(0)=J(\bar s)=0$. The {\em multiplicity} of a conjugate instant is the dimension of the vector space of the Jacobi vector fields vanishing at $0$ and $\bar s$. Two points $p_0$ and $q_0$ on $M$ are said to be {non-conjugate} in $(M,F)$ if $\bar s=1$ is not a conjugate  instant along any geodesic $x\colon [0,1]\to M$ such that $x(0)=p_0$ and $x(1)=q_0$. We observe that on a Randers manifold we can consider other type of Jacobi vector fields, given by variations of geodesics parametrized to have constant speed with respect to the  Riemannian metric, but as commented in Remark \ref{bothjacobi}, they 
generate the same conjugate points  as the classical ones.

The function  $G=F^2$ is smooth outside the zero section but it is only $C^1$ on the whole tangent bundle. It is $C^2$ on $TM$  if and only if it is the square of the norm of a Riemannian metric (see \cite{Warner65}). Hence the lack of regularity on the zero section is a characteristic property of Finsler
metrics.
This has consequences on the level of regularity of the energy functional of a Finsler metric.
It is easy to see that $E$ is a $C^{1,1}$-functional on  $\Omega_{p_0,q_0}(M)$, i. e. it is differentiable with  locally Lipschitz differential (see \cite[Theorem 4.1]{Mercur77}) but it is well known that   $E$ is not a $C^2$-functional on
$\Omega_{p_0,q_0}(M)$. This fact  makes difficult  the application of infinite dimensional methods in Morse theory for Finsler manifolds and indeed approximations  of $\Omega_{p_0,q_0}(M)$ (or of the free loop space $\Omega(M)$, in the closed geodesics problem) by finite dimensional manifolds are commonly used to apply Morse theory to the energy functional  of a Finsler metric, see for instance \cite[Chapter 17]{Shen01} (or \cite{Matthi80} and \cite{BanLon07} in the periodic case).

Nevertheless, in several papers about geodesics on  Finsler manifolds
it is claimed that $E$ is  twice Frechet differentiable in the
$H^1$-topology at any critical point (a geodesic). In particular, in
\cite{MerPal87} the authors prove an extension of the Morse Lemma to
the case of a $C^{1,1}$-functional defined on a Hilbert manifold and
twice Frechet differentiable at any non-degenerate critical point
and then apply their result to cover Morse theory for closed
geodesics of a bumpy Finsler metric on a compact manifold. Moreover
in \cite{MouSou00} the results of \cite{MerPal87} are extended
proving the splitting lemma at a degenerate and isolated critical
point.

Unfortunately, as recently shown  by A. Abbondandolo and M. Schwarz
\cite{AbbSch09}, the energy functional of a Finsler metric is twice
differentiable at a critical point if and only if $F^2(x,y)$ is
Riemannian along the critical point. Although \cite[Proposition
2.3]{AbbSch09} deals with  smooth time-dependent Lagrangians having
at most quadratic growth in the velocities, the argument developed
there works also for the square of a Finsler metric since it does
not require continuity  of the derivatives $\frac{\partial^2
(F^2)}{\partial y^i\partial y^j}$ on the zero section. Moreover the
proof in \cite{AbbSch09} concerns any curve in the manifold
$H^1([0,1],M)$ but it works, with minor modifications, also for
curves satisfying periodic or fixed endpoint boundary conditions.
Without the Morse Lemma, the  computation of the critical groups ,
which are the local homotopic invariants describing  the ``nature''
of an isolated critical point (see \cite{Rothe73}), cannot be
carried out  in an infinite dimensional setting.

In Section~\ref{split} we show (Theorem~\ref{ML}) that, in spite of
the lack of twice differentiability, the splitting lemma holds for
the energy functional $E$ on the infinite dimensional manifold
$\Omega_{p_0,q_0}(M)$.

To this end, we use some ideas of K.-C. Chang who
proved a   splitting lemma (i.e.  an extension of the Morse Lemma to the case of a degenerate critical point) for a $C^2$-functional $J$ defined on a Banach space $X$ immersed continuously as a dense subspace of a Hilbert space $H$ and whose gradient is of the type $\nabla J(x)=x - K\cdot L(x)$, where $L$ is a $C^{1,1}$-map from $H$ to another Banach space $E$ and $K$ is a continuous linear operator from $E$ to $X$ (see \cite{Chang83} and also \cite[Remark 5.1.15]{Chang05}).  Such a result was extended  in \cite{LiLiLi05} for a  $C^{1,1}$-functional on  $H$
which is $C^2$ on an open subset $U$ of  $X$. Similar ideas have also appeared in  M. Struwe's work about the Plateau's problem (cf. \cite{Struwe88})
and in \cite{Ming99}. In particular
the extension of the  splitting lemma to Banach spaces proved in this last paper is suited also for the energy functional of a Finsler metric.

In fact, the energy functional $E$ is $C^2$ on the manifold of the smooth regular curves, having fixed endpoints, endowed with the $C^1$-topology. After a {\em localization procedure} which allows us to work on the Hilbert space $H^1_0([0,1],U)$, $U$ being an open subset of $\R^n$, $n=\mathrm{dim} M$,  the  extension to $H^1_0([0,1],\R^n)$ of the second Frechet derivative of $E$ at a  critical point $\bar x$ is given, with respect to a scalar product $(\cdot,\cdot)$ equivalent to the standard one,  by $A(\bar x)=I+K(\bar x)$, where $I$ is the identity operator  of $H^1_0([0,1],\R^n)$ and $K(\bar x)$ is a  bounded linear operator from $H^1_0([0,1],\R^n)$ to $C^1_0([0,1],\R^n)$. More important, the gradient of $E$ evaluated at the curves in $D\subset C^1_0([0,1],U)$ is  a field in $C^1_0([0,1],\R^n)$. Here $D$ is the open subset of $C^1_0([0,1],U)$ which corresponds to the curves where the localized Lagrangian is regular (see the beginning of Section~\ref{split}).
Using the scalar product $(\cdot,\cdot)$ and the operator $A(\bar x)$ to represent $E$, we obtain a splitting lemma for $E$ restricted to $C^1_0([0,1],U)$ (Theorem~\ref{ML}) and the Morse relations for the  geodesics connecting two non-conjugate points in $(M,F)$.

\end{section}
\begin{section}{The splitting lemma for the energy functional of a Finsler metric and the Morse relations}\label{split}
By using a localization argument (see  \cite[Appendix A]{AbbFig07}), we can assume that the energy functional $E$ is given in a coordinate system  of the manifold $\Omega_{p_0,q_0}(M)$ by
\beq\label{eloc}
\tilde E(x)=\frac 12\inte\tilde{G}(s,x(s),\dot x(s))\de s,
\eeq
where $\tilde G$ is a ``time-dependent'' (non-homogeneous) Lagrangian  defined on an open subset $U$ of $\R^n$, $n=\dim M$.
The localization argument works as follows.
Assume that $\bar x\colon [0,1]\to M$ is a  differentiable curve of the Finsler metric $F$ connecting the points $p_0$ and $q_0$. Let $\exp$ be the exponential map of the auxiliary Riemannian metric $h$, $\mu(p)$ be the injectivity radius of the point $p$ in $(M,h)$
and $\rho=\inf\{\mu(p):p\in \bar x([0,1])\}$.
Let $[0,1]\ni s\rightarrow{\mathbf E}(s)=(E_1(s),\ldots,E_n(s))$ be a parallel orthonormal frame along $\bar x$, $P_s:\R^n\rightarrow T_{\bar x(s)}M$  defined as $P_s(q_1,\ldots,q_r)=q_1 E_1(s)+\ldots +q_n E_n(s)$ and consider the  Euclidean open ball $U$ of radius $\rho/2$ and the map
$\varphi(s,q)={\rm exp}_{\bar x(s)}P_s(q)$.  The map $\varphi_s:U\rightarrow M$, defined as $\varphi_s(q)=\varphi(s,q)$,  is  injective with invertible differential $\de \varphi_s(q)$, for every $s\in[0,1]$ and  $q\in U$.

The Lagrangian $\tilde G\colon[0,1]\times U\times\R^n\to\R$ is defined as
\begin{equation}\label{lagrangian}\tilde G(s,q,y)=F^2(\varphi(s,q),\de \varphi(s,q)[(1,y)]).
\end{equation}
It is continuous on $[0,1]\times U\times\R^n$.
The lack of regularity of $F^2$ on the zero section of $TM$ is inherited by $\tilde G$
on the set $Z\subset[0,1]\times U\times\R^n$ given by all the points $(s,q,y)$ such that $\de \varphi(s,q)[(1,y)]=0$. Observe that for each $(s,q)\in[0,1]\times U$ there is only  one $y\in \R^n$ such that
$\de \varphi(s,q)[(1,y)]=0$. In fact, $\de \varphi(s,q)[(1,y)]=\partial_s \varphi(s,q)+\partial_q\varphi(s,q)[y]$, where $\partial_s \varphi(s,q)$ and $\partial_q \varphi(s,q)$ are the partial differentials of $\varphi$ with respect to the $s$ and $q$ variables; as $\partial_q\varphi(s,q)$ is injective, $y\in \R^n$ is the only vector such  that
$\partial_q\varphi(s,q)[y]=-\partial_s \varphi(s,q)$.

Since $F^2$ is fiberwise strictly convex, we have that
$\tilde G_{yy}(s,q,y)$ is positive definite for all $(s,q,y)\in [0,1]\times U\times\R^n\setminus Z$.

Define the map
\beq\label{varfistar}
\varphi_*\colon H^1_0([0,1],U)\to \Omega_{p_0,q_0}(M),\quad\quad \varphi_*(\xi)(s)=\varphi(s, \xi(s)).\eeq
Hence
\beq\label{tildee}
\tilde E=E\circ\varphi_*.
\eeq
Observe  that the constant function $0\in H^1_0([0,1],U)$ is mapped by $\varphi_*$ to the geodesic $\bar x$.

Let $D$ be the open subset of   $C^1_0([0,1],U)$ containing all the curves $x$ with the property
$(s,x(s),\dot x(s))\not \in Z$, for all $s\in[0,1]$.

By a standard argument, it can be proved that $\tilde E$ is twice Frechet  differentiable at $x$ in $D$ endowed with the $C^1$-topology.
\begin{lemma}\label{gateauxder}
$\tilde{E}$ admits second Frechet derivative  $D^2\tilde E(x)$ at a curve $x\in D$, with respect to the $C^1$-topology
and it is given by
\bal
\lefteqn{D^2\tilde E(x)[\xi_1,\xi_2]=}&\nonumber\\
&=\frac12\inte\big(\tilde G_{qq}(s,x,\dot x)[\xi_1,\xi_2]+\tilde G_{yq}(s,x,\dot x)[\xi_1,\dot \xi_2]\big)\de s\nonumber\\
&\quad+\frac 12\inte\big(\tilde G_{qy}(s,x,\dot x)[\dot \xi_1,\xi_2]\de s+\tilde G_{yy}(s,x,\dot x)[\dot \xi_1,\dot \xi_2]\big)\de s.\label{gateaux}\eal
\end{lemma}
Observe that the right-hand  side of \eqref{gateaux} can be extended to a bounded symmetric bilinear form $B$ on $H^1_0([0,1],\R^n)$.

Observe also that the vector fields $\nu$ in the kernel $N$ of $B$ in $H^1_0([0,1],\R^n)$ correspond to the  Jacobi fields along the geodesic $x$, vanishing at the endpoints. Therefore they are smooth
and $N$ is finite dimensional.

Since $\tilde G$ is fiberwise strictly convex, the bilinear form
\beq\label{scalarproduct}
(\xi_1,\xi_2)\mapsto \frac{1}{2}\inte \tilde G_{yy}(s,0,0)[\dot \xi_1,\dot \xi_2]\de s\eeq
defines a scalar product $(\cdot,\cdot)$ on $H^1_0([0,1],\R^n)$ which is equivalent to the standard one.
\bl\label{A}
Let $B$ be the extension of $D^2\tilde{E}(0)$ to $H^1_0([0,1],\R^n)$. There exists  a bounded linear operator $A\colon H^1_0([0,1],\R^n)\to H^1_0([0,1],\R^n)$ of the type $A=I+K$ where $I$ is the identity operator and $K\colon H^1_0([0,1],\R^n)\to H^1_0([0,1],\R^n)$ is a bounded linear operator,  such that $B$ is represented with respect to the scalar product \eqref{scalarproduct} by $A$.   Moreover the range of $K$ is contained in $ C^1_0([0,1],\R^n)$ and, as an operator $H^1_0([0,1],\R^n)\to C^1_0([0,1],\R^n)$, $K$ is bounded as well.
\el
\begin{proof}
From \eqref{gateaux}, $K$
 is the sum of the bounded linear operators 
\[K_i\colon H^1_0([0,1],\R^n)\to C^1_0([0,1],\R^n), \quad i=1,2,3,\]  defined as follows.
For each $s\in[0,1]$, let  $\tilde G^{yy}(s,0,0)$ be the inverse matrix  of $\tilde G_{yy}(s,0,0)$. For any $\xi\in H^1_0([0,1],\R^n)$ let $K_1\xi$ be the $C^1$-vector field $W_1$ which solves the equation
\[\frac{\de}{\de s}\left(\tilde G_{yy}(s,0,0)\dot W_1\right)=-\tilde G_{qq}(s,0,0)\xi,\]
and vanishes at $s=0,1$, so that
\[\frac12\inte\tilde G_{qq}(s,0,0)[\xi_1,\xi_2]\de s=\frac12\inte\tilde G_{yy}(s,0,0)[\frac{\de}{\de s}\left(K_1(\xi_1)\right),\dot\xi_2]\de s.\]
Hence
\beq\label{w1}\dot W_1=-\tilde G^{yy}(s,0,0)\int_0^s\tilde G_{qq}(\tau,0,0)\xi\de \tau +\tilde G^{yy}(s,0,0) C_1(\xi),\eeq
where $C_1(\xi)$ is the constant vector equal to
\beq\label{c1}
C_1(\xi)=\left(\inte\tilde G^{yy}(s,0,0)\de s\right)^{-1}\left(\inte \tilde G^{yy}(s,0,0)\left(\int_0^s\tilde G_{qq}(\tau,0,0)\xi\de\tau\right)\de s\right)\eeq
(notice that since $\tilde G^{yy}(s,0,0)$ is positive definite for all $s\in [0,1]$, the matrix  $\int_0^1\tilde G^{yy}(s,0,0)\de s$ is positive definite and invertible).

Analogously
$K_2\xi$ and $K_3\xi$ are the curves $W_2$ and $W_3$ in $C^1_0([0,1],\R^n)$ which solve respectively the equations
\bal
&\dot W_2=\tilde G^{yy}(s,0,0)\cdot \tilde G_{yq}(s,0,0)\xi +\tilde G^{yy}(s,0,0)C_2(\xi)\label{w2},\\
&\dot W_3=-\tilde G^{yy}(s,0,0)\int_0^s\tilde G_{qy}(\tau,0,0)\dot \xi\de\tau+\tilde G^{yy}(s,0,0) C_3(\xi)\label{w3},
\eal
where $C_2(\xi)$ is the constant vector equal to
\bal
&\label{c2}C_2(\xi)=-\left(\int_0^1\tilde G^{yy}(s,0,0)\de s\right)^{-1}\left(\inte\tilde G^{yy}(s,0,0)\tilde G_{yq}(s,0,0)\xi\de s\right)\\
\intertext{and}
&\label{c3}C_3(\xi)=\left(\inte\tilde G^{yy}(s,0,0)\de s\right)^{-1}\left(\inte \tilde G^{yy}(s,0,0)\left(\int_0^s\tilde G_{qy}(\tau,0,0)\dot\xi\de\tau\right)\de s\right).
\eal
\end{proof}
\bere
By the compact embedding of $H^1_0([0,1],\R^n)$ in $C^0_0([0,1],\R^n)$,  it follows by \eqref{w1}, \eqref{c1}, \eqref{w2}, \eqref{c2} that $K_1$ and $K_2$ are compact operators in $H^1_0([0,1],\R^n)$, moreover from the Ascoli-Arzel\`a theorem and \eqref{w3}, \eqref{c3}, also $K_3$ is compact and then $A$ is a Fredholm operator in $H^1_0([0,1],\R^n)$ with closed range equal to $N^{\perp}$.
\ere
\bere\label{tildeP}
Since  $N$ is contained in $C^1_0([0,1],\R^n)$, every  $\xi\in C^1_0([0,1],\R^n)$, as a vector field in $H^1_0([0,1],\R^n)$, has projection $P\xi$ on $N^{\perp}$  which is also in   $C^1_0([0,1],\R^n)$. Hence $C^1_0([0,1],\R^n)$ is the topological direct sum of the closed subspaces $N$ and $N^{\perp}\cap C^1_0([0,1],\R^n)$.

Let us define by $\tilde P\colon N\oplus \big(N^{\perp}\cap C^1_0([0,1],\R^n)\big) \to N^{\perp}\cap C^1_0([0,1],\R^n)$ the projection operator.
In the following, we  denote by $\|\cdot\|$ and $\|\cdot\|_{C^1}$ respectively the norm of  $H^1_0([0,1],\R^n)$ endowed with the scalar product \eqref{scalarproduct} and the norm of the $C^1$ topology in $C^1_0([0,1],\R^n)$.
\ere
\bl\label{Arestricted}
The restriction $\tilde A$ of $A$ to the subspace $N^{\perp}\cap C^1_0([0,1],\R^n)$ is an invertible operator $\tilde A\colon  N^{\perp}\cap C^1_0([0,1],\R^n)\to N^{\perp}\cap C^1_0([0,1],\R^n)$ with bounded inverse.
\el
\begin{proof}
Let $\eta\in  N^{\perp}\cap C^1_0([0,1],\R^n)$. Observe that as a curve in $H^1_0([0,1],\R^n)$, $A\eta$ belongs to $N^{\perp}$. Since $A\eta=\eta+K\eta$ and $R(K)\subset C^1_0([0,1],\R^n)$, $A\eta\in C^1_0([0,1],\R^n)$. Moreover
\[\|A\eta\|_{C^1}\leq\|\eta\|_{C^1}+\|K\eta\|_{C^1}\leq \|\eta\|_{C^1}+\|K\|\|\eta\|\leq \|\eta\|_{C^1}+\|K\|\|\eta\|_{C^1}.\]
Therefore $A$ is bounded from $N^{\perp}\cap C^1_0([0,1],\R^n)$ to $N^{\perp}\cap C^1_0([0,1],\R^n)$.
Moreover for any $\bar \eta\in N^{\perp}\cap C^1_0([0,1],\R^n)$ let $\eta\in N^{\perp}$ such that $A\eta=\bar\eta$. Hence $\eta=\bar \eta- K\eta\in C^1_0([0,1],\R^n)$ that is $\tilde A$ is surjective and by the open mapping theorem it has bounded inverse.
\end{proof}
\bl\label{smoothgradient}
Let $x\in D$, then $\nabla\tilde E(x)\in C^1_0([0,1]),\R^n)$. Moreover the map $x\in D\mapsto \nabla\tilde E(x)$ is continuous in the $C^1$ topology.
\el
\begin{proof}
Let $x\in D$, the differential of $\tilde E$ at $x$ in $H^1_0([0,1],U)$ is given by
\[\de \tilde E(x)[\xi]=\frac 12\inte \big(\tilde G_q(s,x,\dot x)\xi+\tilde G_y(s,x,\dot x)\dot \xi\big)\de s\]
for all $\xi\in H^1_0([0,1],\R^n)$.
Recalling that we are using the scalar product \eqref{scalarproduct} on $H^1_0([0,1],\R^n)$,
$\nabla\tilde E(x)$ is the curve $W\in  H^1_0([0,1],\R^n)$ such that $(W,\xi)=\de \tilde E(x)[\xi]$, that is
\[\frac 12\inte \tilde G_{yy}(s,0,0)[\dot W,\dot\xi]=\frac12 \inte\big(\tilde G_q(s,x,\dot x)\xi+\tilde G_y(s,x,\dot x)\dot\xi\big)\de s.\]
Thus
\[\frac 12\inte \tilde G_{yy}(s,0,0)[\dot W,\dot\xi]-\frac12 \inte\left(-\int_0^s\tilde G_q(\tau,x,\dot x)\de \tau+\tilde G_y(s,x,\dot x)\right)\dot\xi\de s=0\]
and this equality is satisfied for all $\xi\in H^1_0([0,1],\R^n)$ if and only if there exists a constant (depending on $x$) $C=C(x)$   such that
\beq\label{dotgrad}\dot W=\tilde G^{yy}(s,0,0)\left(-\int_0^s\tilde G_q(\tau, x,\dot x)\de \tau+\tilde G_y(s,x,\dot x)+C(x)\right).\eeq
As $W$ must vanish at $s=0$ and $s=1$, $C(x)$ has to be equal to
\beq\label{cgrad}
C(x)=\left(\inte \tilde G^{yy}(s,0,0)\de s\right)^{-1}\inte \tilde G^{yy}(s,0,0)\left(\int_0^s\tilde G_q(\tau, x,\dot x)\de \tau-\tilde G_y(s,x,\dot x)\right)\de s.\eeq
From \eqref{dotgrad} and \eqref{cgrad}, we see that $\nabla\tilde E(x)\in C^1_0([0,1],\R^n)$ and using uniform continuity of the vector fields $\tilde G_q(s,q,y)$ and $\tilde G_y(s,q,y)$ we get that
$\|\nabla E(x_n)-\nabla E(x)\|_{C^1}\to 0$ if $x_n\to x$ in the $C^1$ topology.
\end{proof}
We are now ready to prove the splitting lemma for $\tilde E$ at a
critical point. From \eqref{tildee}, since the map $\varphi_*$ is
smooth and injective, we obtain the splitting lemma (or in case the
geodesic $x_0$ is non-degenerate, the Morse lemma) for $E$. A proof   
of the splitting lemma for the energy functional of a Finsler manifold was established in
\cite[Lemma 4.2]{Matthi80}, see also \cite[\S 17.4]{Shen01}, using a
finite dimensional reduction on the manifold of piecewise minimizing
geodesics. Here we present an infinite dimensional proof in the
spirit of the papers of Gromoll and Meyer
\cite{GroMey69a,GroMey69b}, see also \cite{Chang05,LiLiLi05,Ming99}.

\bt\label{ML}
Let $x_0$ be  a geodesic of the Finsler manifold $(M,F)$ connecting two  points $p_0$ and $q_0$ in $M$ and consider the function $\varphi_*$ defined in \eqref{varfistar} associated to $x_0$. Then there exist a ball $B(0,r)$ in $C^1_0([0,1],U)$ centered at $0$, a  local homeomorphism
$\phi\colon B(0,r)\to \phi(B(0,r))\subset D$,  $\phi(0)=0$, a $C^1$ map $h\colon B(0,r)\cap N\to D\cap N^{\perp}$, where $N$ is  the kernel of $A$
such that
\beq\label{splitting}\tilde E(\phi(\xi))=\frac 1 2\big (A\eta,\eta\big)+ \tilde E(\nu+h(\nu)),\eeq
 $\xi=\eta+\nu$ with $\nu\in N$ and $\eta\in N^{\perp}$, where $N^{\perp}$ is the orthogonal of $N$ with respect to $(\cdot,\cdot)$.
\et
\begin{proof}
Consider the equation
\beq\label{key}
\tilde P\cdot \nabla\tilde E(\nu+\eta)=0,\eeq
where  $(\nu,\eta)\in(D\cap N)\times (D\cap N^{\perp})$ and $\tilde P$ was defined in Remark~\ref{tildeP}.
The function $(\nu,\eta)\in(D\cap N)\times (D\cap N^{\perp})\mapsto \tilde F(\nu+\eta)=\tilde P\nabla \tilde E(\nu+\eta)\in N^{\perp}\cap C^1_0([0,1],\R^n)$ is continuous by Lemma~\ref{smoothgradient}, moreover using  \eqref{dotgrad} and \eqref{cgrad}, one can prove by a standard argument that  it is differentiable with respect to the $C^1$-topology with continuous differential. In particular $\tilde F$ is differentiable with respect to $\eta$ and its partial derivative $\de_{\eta} \tilde F(0,0)\colon N^{\perp}\cap C^1_0([0,1],\R^n)\to N^{\perp}\cap C^1_0([0,1],\R^n)$ is the bounded invertible operator $\tilde A$.
Namely, since $\tilde P$ is a bounded linear operator and $\tilde A$ assumes values in $N^{\perp}\cap C^1_0([0,1],R^n)$, it is enough  to prove that
\[\frac{\|\nabla \tilde E(\eta)-\tilde A\eta\|_{C^1}}{\|\eta\|_{C^1}}\To 0,\quad\quad \text{as $\|\eta\|_{C^1}\to 0$}.\]
From \eqref{w1}, \eqref{w2}, \eqref{w3} and \eqref{dotgrad},
we have
\bal
\lefteqn{\frac{\de}{\de s}\left(\nabla\tilde E(\eta)-\tilde A\eta\right)=}&\nonumber\\
&=\tilde G^{yy}(s,0,0)\left(-\int_0^s\tilde G_q(\tau, \eta,\dot \eta)\de \tau+\tilde G_y(s,\eta,\dot \eta)+C(\eta)\right)\nonumber\\
&\quad -\tilde G^{yy}(s,0,0)\tilde G_{yy}(s,0,0)\dot\eta\nonumber\\
&\quad\quad+\tilde G^{yy}(s,0,0)\int_0^s\tilde G_{qq}(\tau,0,0)\eta\de \tau -\tilde G^{yy}(s,0,0) C_1(\eta)\nonumber\\
&\quad\quad\quad-\tilde G^{yy}(s,0,0)\cdot \tilde G_{yq}(s,0,0)\eta -\tilde G^{yy}(s,0,0)C_2(\eta)\nonumber\\
&\quad\quad\quad\quad+\tilde G^{yy}(s,0,0)\int_0^s\tilde G_{qy}(\tau,0,0)\dot \eta\de\tau-\tilde G^{yy}(s,0,0) C_3(\eta)\label{deriv}
\eal
Since the constant curve of constant value $0$ is a critical point of $\tilde E$,  i.e.  $0=\nabla\tilde E(0)$, 
also  the derivative of the curve $\nabla\tilde E(0)$ is constant and equal to zero and then from
\eqref{dotgrad} (with $x=0$) we get
\beq\label{zero}
 0=\tilde G^{yy}(s,0,0)\left(-\int_0^s\tilde G_q(\tau, 0,0)\de \tau+\tilde G_y(s,0, 0)+C(0)\right),\eeq
and we can add the function on the right-hand side of \eqref{zero}   in the equality \eqref{deriv}. By using the mean value theorem, for each $s\in[0,1]$ and to each component of the function
\[t\in[0,1]\mapsto -\int_0^s\tilde G_q(\tau, t\eta(\tau),t\dot \eta(\tau))\de \tau+\tilde G_y(s,t\eta(s),t\dot \eta(s))\]
and the uniform continuity in $[0,1]\times U\times \R^n\setminus Z$ of the second derivatives of the function $\tilde G$, we get that for all $\eps>0$ there exists $\delta>0$ such that for all $\eta\in N^{\perp}\cap C^1_0([0,1],\R^n)$
with $\big\|\eta\|_{C^1}<\delta$
\bmln\Big\|\tilde G^{yy}(s,0,0)\Big(-\int_0^s\tilde G_q(\tau, \eta,\dot \eta)\de \tau+\int_0^s\tilde G_q(\tau, 0,0)\de \tau+\tilde G_y(s,\eta,\dot \eta)\\-\tilde G_y(s,0,0)
-\tilde G_{yy}(s,0,0)\dot\eta+\int_0^s\tilde G_{qq}(s,0,0)\eta\de\tau\\
-\tilde G_{yq}(s,0,0)\eta+\int_0^s\tilde G_{qy}(\tau,0,0)\dot \eta\de\tau\Big)\Big\|_{C^0}<\eps\|\eta\|_{C^1},\emln
where $\|\cdot\|_{C^0}$ is the norm in the $C^0$-topology.
Analogously, since
\[\left(\inte\tilde G^{yy}(s,0,0) \de s\right)^{-1}\inte\tilde G^{yy}(s,0,0)\tilde G_{yy}(s,0,0)\dot\eta\de s=0,\]
recalling \eqref{c1}, \eqref{c2}, \eqref{c3} and \eqref{cgrad}, we have
\bmln
\Big\|\tilde G^{yy}(s,0,0)\big(C(\eta)-C(0)-C_1(\eta)-C_2(\eta)-C_3(\eta)\big)\\
+G^{yy}(s,0,0)\Big(\inte\tilde G^{yy}(\tau,0,0)\de\tau\Big)^{-1}\inte\tilde G^{yy}(\tau,0,0)\tilde G_{yy}(\tau,0,0)\dot\eta\de \tau\Big)\Big\|_{C^0}<\eps\|\eta\|_{C^1}\emln
Therefore, by the implicit function
theorem there exists a $C^1$ map $h\colon B(0,r_1)\cap N\to B(0,\delta_1)\cap N^{\perp}$
such that
\beq\label{proj}
\tilde P\cdot \nabla\tilde E(\nu+h(\nu))=0,\eeq
where $B(0,r_1)$ and $B(0,\delta_1)$ are two balls in $D$ centered at $0$.

Now we consider the Cauchy problem in $H^1_0([0,1],\R^n)$
\beq\label{cauchy}\begin{cases}
\dot \psi(s)=-\dfrac{A\psi(s)}{||A\psi(s)||}\\
\psi(0) =u
  \end{cases}\eeq
where $u\in B(0,r_1)$.
Observe that, as $\|\psi(s)-u\|\leq |s|$, we have $\|\psi(s,u)\|\geq \|u\|-|s|$, thus the flow $\psi$ is well-defined for $|s|<\|u\|$ and $\psi(s,u)\in N^{\perp}$ if $u\in N^{\perp}$.
By Lemma~\ref{Arestricted}, we can solve the above ODE  in the Banach space $N^{\perp}\cap C^1_0([0,1],\R^n)$ (observe that, since the function $u\in N^{\perp}\cap C^1_0([0,1],\R^n)\mapsto (Au,Au)$ is continuous with respect to the $C^1$-topology, the right-hand side of the equation in \eqref{cauchy} is also locally Lipschitz in $N^{\perp}\cap C^1_0([0,1],\R^n)\setminus \{0\}$ with the $C^1$-topology).

Let us call $\zeta$ the flow of \eqref{cauchy} in $N^{\perp}\cap C^1_0([0,1],\R^n)$.
By the uniqueness of the solutions of the Cauchy problem \eqref{cauchy}, 
we have $\zeta(s,u)=\psi(s,u)$ for all $s\in [0,\|u\|)$ and $u\in B(0,r_1)$. Moreover the map $(s,u)\mapsto\psi(s, u)$ is continuous, on the subset of $\R\times B(0,r_1)$ where it is defined, with respect to the product topology of $\R$ and  $B(0,r_1)$ with the $C^1$-topology.
Thus we can adapt the proof of the splitting lemma
in \cite[Theorem 5.1.13]{Chang05} to get the thesis. Namely, consider the functions \baln
&\mathcal F(u,\nu)= \tilde E(\nu+\eta) -\tilde E(\nu+h(\nu)),
&\mathcal F_2(u)=\frac 12 (Au,u),\ealn
where $u=\eta-h(\nu)\in N^{\perp}\cap C^1_0([0,1],\R^n)$. Observe that $\mathcal F(0,\nu)=0$ and $\de_u \mathcal F(0,\nu)=\de_{\eta}\tilde E(\nu+h(\nu))$. Since in the $C^1$ topology, \[\de \tilde{E}(x)[\xi]=\frac 12 \inte \big(\tilde G_q(s,x,\dot x)\xi+\tilde G_y(s,x,\dot x)\dot \xi\big)\de s =(\nabla \tilde E(x),\xi),\]
from  \eqref{proj} we get that $\de_{\eta}\tilde E(\nu+h(\nu))[\xi]=0$ for all $\xi\in N^{\perp}\cap C^1_0([0,1],\R^n)$. Observe also that the second Frechet derivative of $\mathcal F$ at $0$ in $D$ with respect to the variable $u$  is equal to
\[\de_u^2\mathcal F(0,0)=\de^2_{\eta}\tilde E(0,0)\]
As before $\de^2 \tilde E(0)[\xi_1,\xi_2]=(A\xi_1,\xi_2)$ and therefore $\de^2_{\eta}\tilde E(0,0)[\xi_1,\xi_2]=(\tilde A\xi_1,\xi_2)=(A\xi_1,\xi_2)$ for all $\xi_1,\xi_2\in N^{\perp}\cap C^1_0([0,1],\R^n)$.

Since $\mathcal F$ is $C^2$ on $D$ with respect to the $C^1$-topology and $\tilde G$ is $C^2$ on $[0,1]\times U\times \R^n\setminus Z$, taking \eqref{gateaux} into account and using the uniform continuity of the second partial derivatives of $\tilde G$,  we can state that  for all $\eps>0$ there exists a ball $B(0,r_2)\subset D$, with $r_2<r_1$ such that
\bal
|\mathcal F(u,\nu)-\mathcal F_2(u)|&=|\mathcal F(u,\nu)-\mathcal F(0,\nu)-\de_u\mathcal F(0,\nu)[u]-\mathcal F_2(u)|\nonumber\\
&=\left|\inte (1-s)\left(\de^2_u \mathcal F(su,\nu)-\de^2_u\mathcal F(0,0)\right)[u,u]\de s\right|\nonumber\\
&<\eps \|u\|^2,\label{dallalto}
\eal
for all $(u,\nu)\in \left(B(0,r_2)\cap N^{\perp}\right)\times \left(B(0,r_2)\cap N\right)$.
Moreover
\bal\label{dalbasso}
|\mathcal F_2(\psi(t,u))-\mathcal F_2(u)|&=\left|\int_0^t\frac{\de}{\de s} \mathcal F_2(\psi(s,u))\de s\right|=\left|\int_0^t\left(\nabla \mathcal F_2(\psi),\dot \psi\right)\de s\right|\nonumber\\
&=\int_0^{|t|}\|A\psi(s)\|\de s\geq C\int_0^{|t|}\|\psi(s)\|\de s\geq
C\left(\|u\||t|-\frac{t^2}{2}\right),
\eal
where $C$ is a positive constant depending only on the spectral decomposition of $A$ in $N^{\perp}$.
As $\mathcal F_2(\psi(t,u))$ is strictly decreasing in $t$, from \eqref{dallalto} and \eqref{dalbasso} we get that,  if $\eps<\frac{C}{4}$,
\[\mathcal F_2(\psi(-t,u))>\mathcal F(u,\nu)>\mathcal F_2(\psi(t,u)),\]
holds,  for all $t$ such that
\[\|u\|\left(1-\sqrt{1-\frac{2\eps}{C}}\right)\leq |t|<\|u\|\]
and for all $u\in B(0,r_2)$. Therefore, by continuity, there exists a unique $\bar t=\bar t(u,\nu)$,
with
\[|\bar t(u,\nu)|\leq  \|u\|\left(1-\sqrt{1-\frac{2\eps}{C}}\right),\]
such that
\beq\label{fine}\mathcal F_2(\psi(\bar t(u,\nu),u))=\mathcal F(u,\nu).\eeq
By the implicit function theorem, the function $\bar t=\bar t(u,\nu)$ has to be continuous in the $C^1$ topology.  Therefore the map $\phi$ is given by the inverse of the map
$\theta= (u,\nu)\in V \mapsto (\vartheta(u,\nu),\nu)$, where
$\vartheta$ is defined as
\[\vartheta(u,\nu)=\begin{cases}0&\text{if $u=0$}\\
\psi(\bar t(u,\nu),u)&\text{if $u\neq 0$}
     \end{cases}\]
and $V=\theta^{-1}(B(0,r))$ where $B(0,r)\subset\theta (B(0,r_2))$.  Eq. \eqref{splitting} then follows from \eqref{fine}.
\end{proof}
\bere
By the localization argument,  the energy functional of a Finsler metric is treated as the action functional of a  Lagrangian which is smooth outside the closed set $Z$  and  it is strictly convex in the velocities. Therefore the splitting lemma above also holds for the action functional of any smooth  Lagrangian of this type or any such a Lagrangian which is  non-smooth only on  a closed subset of $TM$ which does not intersect the  support of the critical point $x$ and its velocity vector field.
\ere
Theorem~\ref{ML} allows us to compute the critical groups of an isolated critical point as, for instance,  in \cite[Corollary 5.1.18]{Chang05}. In particular we can obtain the Morse relations of geodesics connecting two non-conjugate points in a Finsler manifold (Theorem~\ref{Morseluce}).

Our  reference about  Morse theory for  a $C^{1,1}$-functional defined on an  infinite dimensional manifold is \cite{MawWil89}.  Let $\Omega$ be a complete Hilbert manifold and $f\colon\Omega\to \R$ be a  $C^{1,1}$-functional.  Let us denote by $\mathcal K$ the set of the critical points of $f$. Let  $u\in\mathcal K$  and $U$ be a neighborhood of $u$ such that $\mathcal K\cap U=\{u\}$. For each $n\in\N$, let  us denote  by $C_n(f,u)$ the $n$-th singular homology group of the pair $(f^c\cap U, f^c\cap U\setminus \{u\})$ over the  field $\mathbb K$, where $c=f(u)$ and $f^c=f^{-1}((-\infty,c])$. Let $b,a\in\R$, $b>a$. We denote by $M(r,f^b,f^a)$ the formal series with coefficients in $\N\cup{+\infty}$ defined by $M(r,f^b,f^a)=\sum_{n=0}^{+\infty}M_n(f^b,f^a)r^n$, where
$M_n(f^b,f^a)=\sum_{u\in {\mathcal K}\cap f^{-1}([a,b])}\mathrm{dim}\,C_n(f,u)$. We  denote $M(r,f^b,\emptyset)$ and $M_n(f^b,\emptyset)$ by $M(r,f^b)$ and $M_n(f^b)$.
Assume that:
\begin{itemize}
\item[i)] all the critical points of $f$ are isolated,
\item[ii)]$f$ satisfies the Palais-Smale condition, i. e. any sequence
$\{x_n\}\subset\Omega$ such that $f(x_n)$ is bounded and $\de f(x_n)\to 0$ as $n\to +\infty$ admits a convergent subsequence
\item[iii)] $M_n(f^b,f^a)$ is finite for every $n$ and equal to zero for $n$ large enough,
\end{itemize}
then there exists a polynomial $Q(r)$, with nonnegative integer coefficients,  such that $M(r,f^b,f^a)=P(r,f^b,f^a)+(1+r)Q(r)$ where
$P(r,f^b,f^a)$ is the Poincar\'e polynomial of the pair $(f^b,f^a)$, i. e. $P(r,f^b,f^a)=\sum_{n=0}^{+\infty}B_n(f^b,f^a)r^n$, where
$B_n(f^b,f^a)$  is the dimension of the $n$-th singular homology group of the pair $(f^b,f^a)$ over the field $\mathbb K$.

Observe that under the assumptions $\mathrm{i)}$ and $\mathrm{ii)}$, $f$ has only a finite number of critical points on the strip  $f^{-1}([a,b])$.
If in addition to $\mathrm{i)-iii)}$, we have also that
\begin{itemize}
\item[iv)]$f$ is bounded from below,
\end{itemize}
then, choosing $a<\inf_{\Omega} f$, we get
\[
M(r,f^b)=P(r,f^b)+(1+r)Q(r).
\]
\bt\label{Morsefinsler}
Let $(M,F)$ be a Finsler manifold, and $p_0,q_0$  be two non-conjugate points in $(M,F)$ and assume  that $F$ is forward or  backward complete. Then there exists a formal series $Q(r)$ with coefficients in $\N\cup\{+\infty\}$ such that
\beq\label{relazioni}
\sum_{x\in \Gamma}r^{\mu(x)}=P(r,\Omega_{p_0,q_0}( M))+(1+r)Q(r),
\eeq
where  $\Gamma$ is the set of all the  geodesics connecting $p_0$ to $q_0$ and $\mu(x)$ is the number of conjugate instants, counted with their multiplicity, along the geodesic $x$.
\et
\begin{proof}
Since the points $p_0$ and $q_0$ are non-conjugate in $(M,F)$, any  critical point $x$ of $E$
in $\Omega_{p_0,q_0}(M)$ is isolated and $A$ has zero null space.

If $(M,F)$ is  forward or backward complete then  $E$ satisfies the Palais-Smale condition on $\Omega_{p_0,q_0}(M)$ (see \cite[Theorem 3.1]{CaJaMa07a}) and it is bounded from below.

Using  Theorem~\ref{ML} we can compute the critical group $C_n(E,x)$.
Let $\mathcal O_*$ be the image of the map $\varphi_*$ in \eqref{varfistar}
associated to the critical point $x$ and consider the functional $\tilde E$
in \eqref{eloc} associated to $\varphi_*$. Since the critical point $x$ is non-degenerate, by Theorem~\ref{ML}, there exists
a  local homeomorphism
$\phi\colon B=B(0,r)\to \phi(B)\subset D$ such that  $\phi(0)=0$ and
\[\tilde E(\phi(\xi))=\frac 1 2 \big (A\xi,\xi\big)+ \tilde E(0).\]
Let  $O=\phi(B)$ and consider the deformation $\psi\colon O\times[0,1]\to O$ defined as
\[\psi(\xi,t)=\phi((1-t)\xi_++\xi_-),\]
where $\xi_++\xi_-=\phi^{-1}(\xi)$ and $\xi_+\in H_+$ and $\xi_-\in H_-$, $H_+$ and $H_-$ being the positive and the negative space  of $A$ according to its spectral decomposition in $H^1_0([0,1],\R^n)$ endowed with the scalar product \eqref{scalarproduct}.
Since $A$ is a compact deformation of the identity operator (see the proof of Lemma~\ref{A}), we know that $H_-$ is finite dimensional.

Then $\psi$ is a deformation retract of $\tilde E^c_{|X}\cap O$ to $\tilde E^c_{|X}\cap O_-$, where $O_-=\phi(B\cap H_-)$, $X=C^1_0([0,1],U)$ and $c=\tilde E(0)$.
Therefore
\beq\label{uno} H_n(\tilde E^c_{|X}\cap O,\tilde E^c_{|X}\cap O\setminus\{0\})=H_n(\tilde E^c_{|X}\cap O_-,\tilde E^c_{|X}\cap O_-\setminus\{0\})=\delta_{n,k}\mathbb K,\eeq
where $k$ is the  index of $A$ as a bilinear form on $C^1_0([0,1],\R^n)$ or equivalently on $H^1_0([0,1],\R^n)$, that is $k=\dim (H_-\cap X)=\dim H_-$, and
$\delta_{n,k}$ is the Kronecker's delta.
By  \cite[Theorem 41.1, Theorem 43.2]{Matsum86}, we also have  $k=\mu(x)$.

Since $C^1_0([0,1],\R^n)$ is immersed continuously in $H^1_0([0,1],U)$, by the excision property of the singular relative homology groups  we have
\beq\label{due} H_n(\tilde E^c_{|X}\cap O,\tilde E^c_{|X}\cap O\setminus\{0\})=
 H_n(\tilde E^c_{|X}\cap \tilde O^*,\tilde E^c_{|X}\cap \tilde O^*\setminus\{0\})
\eeq
where $\tilde O^*$ is any neighborhood of $0$ in $H^1_0([0,1],U)$.
On the other hand by \cite[Theorem 16 and Theorem 17]{Palais66a},
and the fact that the map $\varphi_*$ is a homeomorphism, we get (see \cite{CaJaMa12}):
\beq\label{tre}
 H_n(\tilde E^c_{|X}\cap \tilde O^*\!,\tilde E^c_{|X}\cap \tilde O^*\setminus\{0\})= H_n(\tilde E^c\cap \tilde O^*\!,\tilde E^c\cap \tilde O^*\setminus\{0\})=  H_n(E^c\cap  O^*\!, E^c\cap  O^*\setminus\{x\})\eeq
where $O^*=\varphi_*(\tilde O^*)$.
Putting together Eqs. \eqref{uno}--\eqref{tre}, we get
\[C_n(E,x)= \delta_{n,k}\mathbb K.\]

Therefore the assumptions $\mathrm{i)-iv)}$ are satisfied and
the Morse relations
\beq\label{rels} \sum_{z\in \Gamma\cap E^b}r^{\mu(z)}=P(r,E^b)+(1+r)Q(r)\eeq
hold. Finally, arguing as in the proof of \cite[Theorem 1.7]{GiMaPi98},  we can pass to the limit on  both sides of \eqref{rels}  obtaining \eqref{relazioni}, as $b\to +\infty$.
\end{proof}
\begin{remark}\label{palais-smale}
Although the distance  \eqref{finsleriandist}  associated to a Finsler metric is not a true distance due to the lack of symmetry, we can define a symmetric distance  as
\begin{equation}\label{symdist}
\dist_s(p,q)=\frac 12 (\dist(p,q)+\dist(q,p)),
\end{equation}
for every $p,q\in M$.
Let us observe that if the closed balls for the symmetrized distance  are compact, then the energy functional of the Finsler metric satisfies the Palais-Smale condition. This fact came out when studying the relation between causality and completeness of Fermat metrics (see \cite[Theorem 4.3 and Theorem 5.2]{CaJaSa09}). This condition is equivalent to have compact intersection $\bar{B}^+(p,r)\cap\bar{B}^-(p,r)$ for every $p\in M$ and $r>0$, where $\bar{B}^+(p,r)=\{q\in M:\, \dist(p,q)\leq r\}$ and $\bar{B}^-(p,r)=\{q\in M:\, \dist(q,p)\leq r\}$ (see \cite[Proposition 2.2]{CaJaSa09}). If $F$ is forward or backward complete, the Finslerian Hopf-Rinow theorem implies that $F$ satisfies the above condition, but the reciprocal is not true (see  \cite[Example 4.6]{CaJaSa09} for an example of a Finsler metric with compact symmetrized closed balls that is neither forward nor backward complete). In fact, the proof of the Palais-Smale condition for the Finslerian energy functional  written out in \cite[Proposition 3.1]{
CaJaMa07a} works also under the above equivalent conditions (for example $\bar{B}^+(p,r)\cap\bar{B}^-(p,r)$ compact for every $p\in M$ and $r>0$). For further details, see the comments before Theorem 5.2 in \cite{CaJaSa09}.  Then the Morse relations for geodesics connecting two non-conjugate points on a Finsler manifold $(M,F)$ hold also under the more general assumption that the closed balls of the symmetrized distance associated to $F$ are compact.
\end{remark}

\bere
The restriction to $C^1$ curves, whose images are in a neighborhood of a given geodesic $x$, can be performed also for periodic boundary conditions. In the Finsler case, we have to take into account the equivariant action of $SO(2)$ on the free loop space $\Omega(M)$. Then a proof of the Morse relations for closed geodesics  of a Finsler metric might be obtained along the same lines of  \cite[Lemma 4]{GroMey69a},
considering the intersection of a tubular neighborhood of an isolated  critical orbit $SO(2)x$ in $\Omega(M)$ with the Banach manifold $C^1(S,M)$.
\ere
\end{section}
\begin{section}{Morse theory of lightlike geodesics}
A conformally standard stationary spacetime  is a Lorentzian manifold $(\mathcal M,g)$ such that
$\mathcal M=\mathcal M_0\times \R$ and
\beq\label{conf}
g(x,t)[(y,\tau),(y,\tau)]=a(x,t)\big(g_0(x)[y,y]+2g_0(x)[\delta(x),y]\tau -\beta(x)\tau^2),
\eeq
where $(x,t)\in \mathcal M_0\times \R$,  $(y,\tau)\in T_x\mathcal M_0\times\R$, $g_0$ is a Riemannian metric on $\mathcal M_0$ and $\delta$, $\beta$ and $a$ are, respectively, a smooth vector field on $\mathcal M_0$,  a smooth positive function on $\mathcal M_0$ and a smooth positive function on $\mathcal M$.
We  denote by $\tilde{g}_0$ the conformal Riemannian metric $g_0/\beta$.
Since lightlike geodesics and conjugate points along lightlike geodesics are preserved under conformal changes of the metric we can  divide \eqref{conf}  by $a\beta$, and so we can assume that the metric $g$ is given by
\begin{equation}\label{l}
g(x,t)[(y,\tau),(y,\tau)]=\tilde{g}_0(x)[y,y]+2\tilde{g}_0(x)[\delta(x),y]\tau -\tau^2.
\end{equation}
By definition, a  smooth lightlike curve $[a,b]\ni s\to\gamma(s)=(x(s),t(s))\in\mathcal M$ has to satisfy the
equation
\[
 \tilde{g}_0(x)[\dot x,\dot x]+2\tilde{g}_0(x)[\delta(x),\dot x]\dot t -\dot t^2=0,
\]
and therefore the derivative of the $t$ component is given by
\beq\label{fpnull}
\dot t=\sqrt{\tilde{g}_0(x)[\dot x,\dot x]+\tilde{g}_0(x)[\delta(x),\dot x]^2}+\tilde{g}_0(x)[\delta(x),\dot x]\eeq
or
\[\dot t=-\sqrt{\tilde{g}_0(x)[\dot x,\dot x]+\tilde{g}_0(x)[\delta(x),\dot x]^2}+\tilde{g}_0(x)[\delta(x),\dot x].\]
Notice that in the first case $\dot t>0$ (the lightlike curve is future-pointing) while in the second one $\dot t<0$ (the lightlike curve is past-pointing).
The right-hand side of the first equation and minus the right-hand side of the second one  define two Randers  metrics on $\mathcal M_0$ that  are denoted, respectively, by $F$ and $F_-$:
\begin{align}\label{Fermatmetric}
F(x,y)&=\sqrt{\tilde{g}_0(x)[y,y]+\tilde{g}_0(x)[\delta(x),y]^2}+\tilde{g}_0(x)[\delta(x),y]\\\nonumber
F_-(x,y)&=\sqrt{\tilde{g}_0(x)[y,y]+\tilde{g}_0(x)[\delta(x),y]^2}-\tilde{g}_0(x)[\delta(x),y]
\end{align}
Such Randers metrics play an important role in the study of lightlike and timelike ge\-o\-des\-ics on a conformally standard stationary spacetime as we  see below and, moreover, they give a lot of information about the causal structure of such type of spacetimes (see \cite{ CaJaMa07a}). As in \cite{CaJaMa07a, CaJaSa09} we  call the Randers metric $F$ the {\em Fermat metric} associated to $(\mathcal M,g)$.

Lightlike geodesics connecting an event $p\in \mathcal M$ with a timelike curve $\gamma \colon (a,b) \to \mathcal M$ can be characterized by a variational principle (a {\em Fermat Principle}) stating that, among all the future-pointing (or past-pointing) lightlike curves $z\colon[0,1]\to\mathcal M$ such that $z(0)=0$ and $z(1)\in \gamma((a,b))$, the lightlike geodesics are all and only the curves making stationary the arrival time functional $T$ that is  the functional $z\mapsto T(z)=\gamma^{-1}(z(1))$. This is a fairly  well known fact since the beginning of general relativity, but a precise formulation with the above generality and a rigorous proof was given only in the '90s by I. Kovner and V. Perlick (see \cite{Kovner90,Perlic90}).
In the case of a conformally standard stationary spacetime,  if we consider an observer whose world line $\gamma$ is the vertical line $\R\ni s\to (x_1,s)\in{\mathcal M}$, the arrival time $T$ coincides with the value of the global time coordinate $t$ at  the endpoint of the curve $[0,1]\ni s\to z(s)=(x(s),t(s))\in {\mathcal M}$. Therefore, for a future-pointing lightlike curve, from \eqref{fpnull} we get
\[T(z)\equiv T(x)=t_0+\inte F(x,\dot x)\de s,\]
hence $T(z)$ is equal, up to an additive constant, to the length with respect to  $F$ of the projection of $z$ on $\mathcal M_0$. The Kovner's Fermat principle can be formulated as follows (see \cite[Theorem 4.4]{CaJaMa07a})
\bpr[Fermat's principle]\label{fermatp}
Let $(\mathcal M,g)$ be a  standard stationary spacetime,
$p=(p_0,t_0)\in\mathcal M$,  $[0,1]\ni s\to\gamma(s)=(q_0,s)\in \mathcal M$, $p_0, q_0\in\mathcal M_0$.  A curve $[0,1]\ni s\to z(s)=(x(s),t(s))\in\mathcal M$ is a future-pointing lightlike geodesic  of $({\mathcal M},g)$ as in \eqref{l} if and only if $[0,1]\ni s\to x(s)\in \mathcal M_0$ is a geodesic of the Fermat metric $F$, parametrized  with constant Riemannian speed $\tilde{g}_0(x)[\delta(x),\dot x]^2+\tilde{g}_0(x)[\dot x,\dot x]=\mathrm{const.}$, and $t(s)$ is given by
\[t(s)=t_0+\!\int_0^sF(x,\dot x)\de v.\]
\epr
By the Fermat's principle the search of lightlike geodesics in a stationary spacetime can be reduced to the search of geodesics in the Finsler manifold $(\mathcal M_0,F)$.

Let $z=(x,t)\colon[0,1]\to\mathcal M$ be a future-pointing lightlike geodesic. By  Proposition~\ref{fermatp}, $x$ is a geodesic in $(\mathcal M_0,F)$. We  denote by $\mu(z)$ (resp. $\mu(x)$)  the {\em geometric index} of $z$  (resp. $x$), that is the number of conjugate points  along $z$ (resp. $x$) counted with their multiplicity.

We recall that on a Lorentzian manifold $(\mathcal M,g)$ the notions of Jacobi vector  field, conjugate instant and non-conjugate points are given, as on a Riemannian manifold, using the Levi-Civita connection and the  Riemannian curvature tensor (see for instance  \cite{One83}).

We are going to  show that the geometric index of $z$  coincides with the geometric index of its spatial projection $x$ as a geodesic of the Fermat metric. This fact allows us to bring the Morse theory for Finsler geodesics to the Morse theory of lightlike geodesics.
\bt\label{index}
Let $(\mathcal M,g)$ be a conformally standard stationary spacetime, $[0,1]\ni s\to z(s)=(x(s),t(s))\in {\mathcal M}$ be a future-pointing lightlike geodesic. Let $F$ be the Fermat metric associated to $(\mathcal M,g)$.
Then
the points $x(0)$ and $x(1)$ are non-conjugate along the geodesic $x$ in $(\mathcal M_0,F)$ if and only if the points $z(0)$ and $z(1)$ are non-conjugate along the lightlike geodesic $z$ in
$(\mathcal M,g)$.  Moreover
\beq \mu(z)=\mu(x).\label{indiciuguali}\eeq
\et
\begin{proof}
As conjugate points  of lightlike geodesics are preserved by conformal changes  with their multiplicity, we can consider the metric $g$ as in \eqref{l},  which  can be expressed as
\beq\label{kk}g(x)[v,v]=\alpha (x)[v,v]-(\tau-\alpha(x)[v,\eta])^2,\eeq
where $\alpha(x)[v,v]=\tilde{g}_0(x)[v,v]+\tilde{g}_0(x)[v,\delta(x)]^2$ and $\alpha(x)[v,\eta(x)]=\tilde{g}_0(x)[v,\delta(x)]$ for every $v\in T_xM$. Let $\bar\nabla$ be the Levi-Civita connection of the metric $\alpha$ and consider the $(1,1)$-tensor field $\Omega$ on $\mathcal M_0$ defined as
\[\Omega[\dot
x]=(\bar\nabla\eta)[\dot x]-(\bar \nabla\eta)^*[\dot x],\]
where $(\bar\nabla
\eta)[\dot x]=\bar\nabla_{\dot x}\eta$ and $(\bar\nabla\eta)^*$ is the
adjoint with respect to $\alpha$ of $\bar\nabla\eta$.

Since $\partial_t$ is a Killing vector field for $(\mathcal M,g)$, we know that  for any geodesic $z=(x,t)$ in $(\mathcal M,g)$ there exists a constant $C_z$ such that
\beq \label{Cz} C_z=\dot t-\alpha(x)[\dot x,\eta].\eeq
Then considering variation vector fields having vanishing $t$ component, from \eqref{kk}, one can easily see that the $x$ component of a geodesic $z$ of $(\mathcal M,g)$, as a critical point of the functional \eqref{nrglorentzian}, has  to satisfy the equation
\beq\label{eul-lag}
\bar\nabla_{\dot x}\dot x=-C_z \Omega[\dot x].
\eeq
The linearized equations of this system \eqref{Cz}--\eqref{eul-lag} are
\begin{align}\label{jacobieq}
J''&=-R(J,\dot x)\dot x-C_{J,W}\Omega[\dot x]-C_z(\bar\nabla_J\Omega)[\dot x]-C_z\Omega[J'],\nonumber\\
W'&=C_{J,W}+\alpha(x)[J',\eta]+\alpha(x)[\dot x,\bar\nabla_J\eta].
\end{align}
On the other hand, the Fermat metric can be expressed as $F(x,v)=\sqrt{\alpha(x)[v,v]}+\alpha (x)[v,\eta]$ and its geodesics with constant $\alpha$-Riemannian speed are determined by
\beq\label{randersgeo}
\bar\nabla_{\dot x}\dot x=-C_x \Omega[\dot x],
\eeq
where $C_x=\sqrt{\alpha(\dot x,\dot x)}$ (see \cite[boxed formula at p. 297]{BaChSh00}). The linearized equation of \eqref{randersgeo} is
\begin{align}\label{jacobieq2}
J''&=-R(J,\dot x)\dot x-\frac{\alpha(x(0))[\dot x(0),J'(0)]}{C_x}\Omega[\dot x]-C_x(\bar\nabla_J\Omega)[\dot x]-C_x\Omega[J'].
\end{align}
If $(J,W)$ is a Jacobi field of $z$ satisfying \eqref{jacobieq} with $(J(0),W(0))=(J(s_0),W(s_0))=(0,0)$, then from \eqref{jacobieq},  using  integration by parts,  Eq. \eqref{eul-lag} and the fact that $\alpha(x)[\dot x, J']$ is constant along $x$ (as can be verified by a direct computation, taking into account that the operators $\Omega$ and  $\bar\nabla_J\Omega$ are skew-symmetric), we obtain the following chain of identities:
\begin{align*}
W(s_0)&=s_0C_{J,W}+\int_0^{s_0}\big(\alpha(x)[J',\eta]+\alpha(x)[\dot x,\bar\nabla_J\eta]\big)\de s\\
&=s_0C_{J,W}+\int_0^{s_0}\alpha(x)[\Omega[\dot x],J]\de s\\
&=s_0C_{J,W}-\frac{1}{C_z}\int_0^{s_0}\alpha(x)[\bar\nabla_{\dot x}{\dot x},J]\de s\\
&=s_0C_{J,W}+\frac{s_0}{C_z}\alpha(x(0))[\dot x(0),J'(0)].
\end{align*}
We observe that $C_z\not=0$ because $z$ is lightlike. As $W(s_0)=0$, last formula implies that $C_{J,W}=-\big(\alpha(x(0))[\dot x(0),J'(0)]\big)/C_z$ and therefore $J$ satisfies \eqref{jacobieq2} taking $C_x=C_z$. Analogously, we can show that if $J$ satisfies \eqref{jacobieq2} and has a conjugate instant $s_0$, then we can construct a Jacobi  vector  field $(J,W)$ satisfying the system \eqref{jacobieq} with $C_z=C_x$, $C_{J,W}=-\big(\alpha(x(0))[\dot x(0),J'(0)]\big)/C_z$ and having a conjugate instant in $s_0$. In conclusion, there is a bijection between the  Jacobi vector  fields of $z$ vanishing in $0$ and $s_0$ and the Jacobi vector  fields of $x$ (as a Fermat geodesic) that are zero in $0$ and $s_0$. This concludes the proof.
\end{proof}
\begin{remark}\label{bothjacobi}
The conjugate points of a Fermat geodesic $x$, when it is pa\-ram\-e\-triz\-ed  with constant $\alpha$-Riemannian speed,  coincide with the conjugate points  when it  is pa\-ram\-e\-triz\-ed with constant Fermat speed. Indeed, let $J$ be a Jacobi vector  field and $\Gamma:[0,1]\times (-\varepsilon,\varepsilon)\to {\mathcal M_0}$ be a variation, by means of geodesics pa\-ram\-e\-triz\-ed with constant Finsler speed, generating $J$. Then we can consider a variation of geodesics with constant $\alpha$-Riemannian speed as
$[0,1]\times (-\varepsilon,\varepsilon)\ni (s,w)\to \tilde{\Gamma}(s,w)=\Gamma(\psi_w(s),w)\in {\mathcal M_0}$, where $[0,1]\ni s\to\psi_w(s)\in[0,1]$ is the reparametrization giving geodesics with constant $\alpha$-Riemannian speed. Consequently, the Jacobi vector  field $\tilde{J}$ corresponding to the variation $\tilde{\Gamma}$ can be expressed as $\tilde{J}(s)=\sigma(s)\dot x(s)+J(\psi_0(s))$ for every $s\in[0,1]$ (here $x(s)=\tilde{\Gamma}(s,0)$). We observe that if $[0,1]\ni s\to \sigma(s)\dot x(s)$  is a Jacobi vector field along $x$, then $\sigma$ must be an affine  function. This can be easily seen, using that $\alpha(x)[\dot x, J']$ is constant. It follows that if we take $J$ such that $J(0)=J(t_0)=0$ with $t_0=\psi_0(s_0)$, then $\bar{J}(s)=\tilde{J}(s)-\frac{\sigma(s_0)}{s_0}s\, \dot x(s)$ satisfies  $\bar{J}(0)=\bar J(s_0)=0$ and it is the unique  Jacobi vector field of the type $\bar{J}(s)=\tilde{J}(s)+\sigma(s) \dot x(s)$ satisfying this property. Conversely, if  $\bar J$ is a Jacobi  vector 
field generated by a variation of geodesics with constant $\alpha$-Riemannian speed and such that $\bar J(0)=\bar J(s_0)=0$, there exists a function $\sigma\colon[0,1]\to\R$ and a family of reparametrizations $[0,1]\ni t\to\phi_w(t)\in[0,1]$ such that $J(t) =\sigma(t)\dot x(t)+\bar J(\phi_w(t))$ is a Jacobi  vector  field corresponding to a variation of geodesics parametrized with constant Finslerian speed. Using again the fact that $\sigma(t)\dot x(t)$ is a Jacobi  vector  field if and only if $\sigma(t)$ is an affine function ( which  can be seen now directly by the Jacobi equation in Finsler geometry, see for instance \cite[formula (6.1)]{Shen01}) we conclude as before   that there exists a unique Jacobi field $J$ corresponding to $\bar J$ such that $J(0)=J(t_0)=0$, where $t_0=\phi_0(s_0)$. Therefore there is a bijection between the conjugate points preserving the points in the geodesic and the order of conjugacy.
\end{remark}

We pass now to study  the Morse relations for lightlike geodesics connecting $p=(p_0,0)$ to the curve $\R\ni s\to\gamma(s)=(q_0,s)\in {\mathcal M}_0\times\R$, $p_0, q_0\in \mathcal M_0$.

\bt\label{Morseluce}
Let $(\mathcal M,g)$ be a globally hyperbolic  conformally standard stationary spacetime, $p=(p_0,t_0)\in\mathcal M$ and $\R\ni s\to\gamma(s)=(q_0,s)\in\mathcal M$. Assume that for each $s\in\R$ the points $p$ and $(q_0,s)$ are non-conjugate along every future-pointing lightlike geodesic connecting them. 
Then there exists a formal series $Q(r)$ with coefficients in $\N\cup\{+\infty\}$ such that
\beq\label{relazioniluce}
\sum_{z\in G_{p,\gamma}}r^{\mu(z)}=P(r,\Omega_{p_0,q_0}(\mathcal M_0))+(1+r)Q(r),
\eeq
where  $G_{p,\gamma}$ is the set of all the future-pointing lightlike geodesics connecting $p$ to $\gamma$.
\et
\begin{proof}
Let $F$ be the Fermat metric in \eqref{Fermatmetric}.  From  Proposition~\ref{fermatp}, any geodesic $x$ in  $(\mathcal M_0, F)$ connecting $p_0$ to $q_0$ corresponds to a future-pointing lightlike geodesic $[0,1]\ni s\to z(s)=(x(s),t(s))\in\mathcal M$ connecting $p$ to $\gamma$ and vice versa. From  Theorem~\ref{index},   the points $p_0$ and $q_0$ are non-conjugate in $(\mathcal M_0,F)$ and $\mu(x)=\mu(z)$. Moreover, by \cite[Theorem 4.3]{CaJaSa09} as $(\mathcal M,g)$ is globally hyperbolic, the Fermat metric $F$ has compact symmetrized closed balls.  Then  \eqref{relazioniluce} comes directly from \eqref{relazioni}  and Remark~\ref{palais-smale}.
\end{proof}
\bere
Observe that taking $r=1$ in \eqref{relazioniluce} gives
\[\sum_{n=0}^{\infty}N_n=\sum_{n=0}^{\infty}B_n(\Omega_{p_0,q_0}(\mathcal M_0))+2Q(1),\]
where $N_n$ is the number of future-pointing  lightlike geodesics having index $n$.
If $\mathcal M_0$ is contractible then $B_0(\Omega_{p_0,q_0}(\mathcal M_0))=1$ and $B_n(\Omega_{p_0,q_0}(\mathcal M_0))=0$ for all $n\geq 1$.
Therefore the number of future-pointing lightlike geodesics joining $p$ to $\gamma$ is infinite or odd.
\ere
\bere
The Morse relations of lightlike geodesics connecting $p$ to  $\gamma$ in a standard stationary spacetime were obtained in \cite{FoGiMa95} by using the  functional
\[\tilde  J (x)=\inte g_0(x)[\delta(x),\dot x]\de s+\left(
\inte \big(g_0(x)[\delta(x),\dot x]^2+g_0(x)[\dot x,\dot x]\big)\de s\right)^{\frac{1}{2}},
\]
and the following Fermat principle: a curve $[0,1]\ni s\to z(s)=(x(s),t(s))\in \mathcal M$ is a future-pointing  lightlike geodesic connecting $p=(p_0,t_0)$ and $\R\ni s\to \gamma(s)=(q_0,s)\in\mathcal M $ if and only if $x$ is a critical point of $\tilde J$
and $t(s)=t_0+\int_0^sF(x,\dot x)\de \nu$.  In \cite{FoGiMa95}, it   was also claimed that the  the Morse index  of  a critical point $x$ of $\tilde J$ is equal  to the  geometrical index of the corresponding  lightlike geodesic $z$, but  there is  a gap in the proof of that statement.
\ere
\end{section}
\begin{section}{Morse theory of timelike geodesics}
The reduction of Morse theory of lightlike geodesics connecting a point with a timelike line on a stationary spacetime $(\mathcal M_0\times \R,g)$ to Morse theory of geodesics of a Finsler metric on $\mathcal M_0$ can be also carried out for timelike geodesics. Namely  timelike geodesics can be viewed as projections on $\mathcal M$ of lightlike geodesics in a one-dimensional higher stationary spacetime as follows.

Let $(\mathcal M, g)$ be a standard stationary spacetime (that is $g$ is given by \eqref{conf} and $a(x,t)=1$; since timelike geodesics are not invariant under conformal changes of the metric, this time we cannot divide $g$ by $\beta$). We seek for  timelike geodesics $z\colon [0,\bar s] \rightarrow{\mathcal M}$ connecting a point $(p_0,t_0)\in \mathcal M$ with a timelike curve
$\R\ni s\to\gamma(s)=(q_0,s)\in \mathcal M$ and parametrized with respect to  proper time i.e.  $E_z = g(z)[\dot z,\dot z]=-1$, for all $s\in[0,\bar s]$.

We extend the Riemannian manifold ${\mathcal M}_0$ to the manifold ${\mathcal N}_0={\mathcal M}_0\times\R$ endowed with the metric
$n_0=g_0+\de u^2$ where  $u$ is the natural coordinate on $\R$,  and we associate
to the manifold ${\mathcal N}_0$ the one dimensional higher  Lorentzian manifold  $({\mathcal  N}, n)$, with the metric $n$ defined as
\begin{equation}\label{extend} n(x,u,t)[(y,v,\tau),(y,v,\tau)]=g_0(x)[y,y]+v^2 + 2g_0(x)[\delta(x),y]\tau -\beta(x)\tau^2.
\end{equation}
 Since  $\partial_u$ is a Killing vector field for the metric  $n$, geodesics $\varsigma=(x,u,t)$ in $({\mathcal N},n)$
have to satisfy also the conservation law $n[\dot \varsigma,\partial_u]=\mathrm{const.}$, which  implies that the $u$ component of  a  geodesic
is an affine function.
Moreover the projection $[a,b]\ni s\to z(s) = (x(s),t(s))\in{\mathcal M}$ of  $\varsigma$ is a geodesic for $({\mathcal M},g)$.
In particular lightlike geodesics of the metric $n$ satisfy the following equation
\[g_0[\dot x,\dot x]+ 2g_0[\delta,\dot x]\dot t -\beta\dot t^2=-\dot u^2=\mathrm{const.}\]
Thus  in order to find timelike geodesics $z=(x,t)$ in $({\mathcal M},g)$, parametrized with respect to  proper time,  it is enough to find
lightlike geodesics in $({\mathcal N},n)$ whose $u$ component has derivative equal to $1$.
The Fermat's principle  can be restated in
$({\mathcal N},n)$, reducing future-pointing lightlike geodesics on $({\mathcal N},n)$ to geodesics for the Fermat  metric $\tilde F$ on the manifold ${\mathcal N}_0$, where $\tilde F$ is given by
\[
\tilde F((x,u),(y,v))=\sqrt{\frac{1}{\beta(x)}(g_0[y,y]+v^2)+\frac{1}{\beta(x)^2}g_0[\delta(x),y]^2 }+\frac{1}{\beta(x)}g_0[\delta(x),y] ,
\]
for all $((x,u),(y,v))\in T{\mathcal N}_0$.
Summing up Theorem~\ref{index} and Theorem~\ref{Morseluce}  we get:
\bt\label{timelike}
Let $(\mathcal M,g)$ be a standard stationary spacetime and $[0,\bar s] \ni s \to z(s) = (x(s),t(s)) \in \mathcal M$ be a future-pointing timelike  geodesic  connecting the point $p=(p_0,t_0)$ to the curve $\R\ni s\to\gamma(s)=(q_0,s)\in\mathcal M$, $q_0\in \mathcal M_0$. Let 
$\tilde F$  be the Fermat metric associated to $(\mathcal N, n)$.
Then
the points $(p_0,0)$ and $(q_0,\bar s)$ are non-conjugate along the geodesic $[0,\bar s]\ni s\to\tilde x(s)=(x(s),s)$ in $(\mathcal N_0,\tilde F)$ if and only if the points $p$ and $(q_0,t(\bar s))$ are non-conjugate along the timelike  geodesic $z$ in $(\mathcal M,g)$. Moreover
\[\mu(z)=\mu(\tilde x).\]
If  for each $s\in\R$, the points $p$ and $(q_0,s)$ are non-conjugate along every future-pointing timelike geodesic, parametrized with respect to proper time on the interval $[0,\bar s]$ and connecting them, and $(\mathcal M,g)$ is globally hyperbolic, then there exists a formal series $Q(r)$ with coefficients in $\N\cup\{+\infty\}$ such that
\[\sum_{z\in \mathcal T_{p,\gamma}}r^{\mu(z)}=P(r,\Omega_{p_0,q_0}(\mathcal M_0))+(1+r)Q(r),
\]
where  $\mathcal T_{p,\gamma}$ is the set of the future-pointing timelike geodesics $[0,\bar s] \ni s \to z ( s ) = ( x ( s ), t ( s ) )\in \mathcal M$ parametrized with respect to proper time  and such that  $z(0)=p$ and $x(\bar s)=q_0$.
\et
\begin{proof}
The first part of the theorem comes arguing as in Theorem~\ref{index}, observing that a Jacobi vector field $\xi=(U,\Upsilon)$ along the  lightlike geodesic $[0,\bar s]\ni s\to\varsigma(s)=(x(s),t(s),s)$ in $(\mathcal N,n)$, with vanishing endpoints, has $\Upsilon$ component equal to $0$ and $U$ component which is a Jacobi vector field along the timelike geodesic $z$.

The second part comes arguing as in Theorem~\ref{Morseluce}, after having observed that  if  the Fermat metric $F$ associated to $(\mathcal M,g)$ has compact symmetrized closed balls the same holds for $\tilde F$. Namely if $\{(x_n,u_n)\}\subset {\mathcal N}_0$ is contained in a symmetrized closed ball of center $(x,u)\in\mathcal M_0\times\R$ and radius $r>0$, then it is easy to see that $x_n$ is contained in the symmetrized closed ball for $F$ of center $x$ and radius $r$, which is compact because $(\mathcal M,g)$ is globally hyperbolic (see \cite[Theorem 4.3 ]{CaJaSa09}). Therefore, there is a subsequence $x_{n_k}$ of $x_n$ that converges
 and $\beta$ is bounded on this subsequence. Thus also $u_{n_k}$ admits a convergent subsequence $u_{n_l}$ in \R and therefore
$\{(x_{n_l},u_{n_l})\}$ converges.
Finally, observe that the manifold  $\Omega_{(p_0,0),(q_0,\bar s)}(\mathcal N_0)$ is homotopically equivalent to $\Omega_{p_0,q_0}(\mathcal M_0)$.
\end{proof}
\end{section}
\begin{section}{Appendix A}
In a  paper about Morse theory of causal geodesics in a globally hyperbolic Lorentzian manifold \cite{Uhlenb75}, K. Uhlenbeck introduced the following functional
\beq\label{tau}
J(z)=\inte\big(g(z)[\dot z,\dot z]+\big(\tfrac{\de Pz}{\de s}\big)^2\big)\de s,\eeq
defined on the set of piecewise differentiable  curves on $\mathcal M$ satisfying the constraint $g(z)[\dot z,\dot z]=0$  and the boundary conditions $z(0)=p$, $z(1)=(q_0,Pz(1))$, where $P\colon \mathcal    M\to \R$ is the natural projection on $\R$, and  proved that critical points of $J$ are all and only the lightlike geodesics connecting $p$ to the line $s\mapsto(q_0,s)$.

In this appendix we study, whenever $\mathcal M$ is conformally standard stationary, the relation between $J$ and the energy functional $E$ of the Fermat metric. We show (Proposition~\ref{morse}) that  the Morse index $m_J(z)$ of a critical point $z=(x,t)$ of $J$ is equal to
the Morse index $m_E(x)$ of $x$ as a critical point of $E$.

This fact provides a variational and alternative proof of the equality \eqref{indiciuguali} since
by Theorem~\ref{Uhl}, $\mu(z)= m_J(z)$
and by  \cite[Theorem 41.1 and Theorem 43.2]{Matsum86} we also have  $\mu(x)=m_E(x)$.

As K. Uhlenbeck observed, the constraint equation $g(z)[\dot z,\dot z]=0$ does not  define a smooth submanifold of the set of piecewise differentiable  curves  in $\mathcal M$, however $J$ is differentiable if viewed as a functional on the set of piecewise differentiable regular (i.e $\dot z(s)\neq 0$, where it is defined) curves on $\mathcal M_0$. This can be done after solving  the constraint equation with respect to $\dot t$, for any fixed $x$. The Lorentzian metric  considered in \cite{Uhlenb75} is  of the type $g(x,t)[(y,\tau),(y,\tau)]=g_0(x,t)[y,y]-\tau^2$ where,
for any $t\in\R$,  $g_0(\cdot,t)$ is a Riemannian metric on $\mathcal M_0$. The solutions $t_x$ of the differential equation
\[\dot t=\sqrt{g_0(x,t)[\dot x,\dot x]},\]
arising from the constraint equation, are defined on the whole interval $[0,1]$ if a rather  technical growth assumption on the    metric $g_0$ is  fulfilled. The critical points of $J$ are exactly the lightlike geodesics $[0,1]\ni s\to z(s)=(x(s),t(s))\in\mathcal M$ connecting the point $p$ to the curve $\R\ni s\to\gamma(s)=(q_0,s)\in\mathcal M$, parameterized with $\dot t$ constant. Moreover she proved the following Morse index theorem
(see \cite[Lemma 4.2]{Uhlenb75}):
\bt[Uhlenbeck]\label{Uhl}
$J$ is twice Gateaux differentiable at any critical point (a lightlike geodesic).
Its second derivative at a critical point $z$ is given by
\bml
D^2J(z)[U,V]=
\inte\big(\tfrac{\de}{\de s}( g(z)[\nabla P,U])\tfrac{\de}{\de s}( g(z)[\nabla P,V])\\
+c(s)\big(g(z)[\nabla_{\dot z} U,\nabla_{\dot z}V]-g(z)[R(\dot z, U)\dot z,V]\big)\big)\de s,
\label{hessJ}\eml
where $U$ and $V$ are piecewise smooth vector fields along $z$ such that $U(0)=V(0)=0=V(1)=U(1)$, $g(z)[\dot z,U]= g(z)[\dot z,V]=0$, $c$ is a function such that $\zeta=c(s)\dot z$, $\zeta$  the parallel transport of $\dot z(1)$ along $z$, $\nabla P$ is the gradient of $P$ and $R$ is the curvature tensor of $(\mathcal M,g)$.
A critical point is non degenerate if and only if its endpoints are non-conjugate. The index of a critical point is equal to its  geometrical index as a lightlike geodesic, that is the number of conjugate points counted with their multiplicity.
\et
\bere
The above theorem is based on the existence of a global time function $P$  and on the fact that $z$ is a lightlike geodesic.
It does not depend  on the  form of the metric $g$, neither on assumptions on the metric coefficients, nor on topological assumptions as global hyperbolicity.
This fact was exploited in the paper \cite{GiMaPi98}.
The key point in the proof of Theorem~\ref{Uhl} is that, if $z$ is a lightlike geodesic, the bilinear form at the right-hand side of \eqref{hessJ} is a compact perturbation of a positive definite invertible operator on the $H^1$ completion of the space of piecewise smooth vector fields $U$ along $z$ satisfying the condition $g[\dot z, U]=0$.

Thus the Morse index of $J$ at a critical point is finite and is equal to the  sum of the dimensions of the kernels of the above bilinear forms along $z(s)$ which  are isomorphic to the space of Jacobi  vector  fields along  $z$ which vanish at the initial point and in some other point  $\bar s\in (0,1)$.
\ere

In the following proposition  the equality between the Morse index of $J$   and $E$ is stated.
\bpr\label{morse}
Under the assumptions of Theorem~\ref{index}, we have that
\beq\label{Morseindex}m_J(z)=m_E(x),\eeq
where  $m_J(z)$ and $m_E(x)$ are respectively the Morse indexes of the functionals $J$ and $E$ at their critical points $z$ and $x$.
\epr
\begin{proof}
Let us denote by $I$ the functional given by \eqref{tau} defined on the manifold $\Omega_{p,\gamma}(\mathcal M)$ of $H^1$- curves in $\mathcal M$ connecting $p$ to $\gamma(\R)$.
Observe that $J$ is equal to the functional $I$ restricted to the set $\Lambda_{p,\gamma}(\mathcal M)\subset\Omega_{p,\gamma}(\mathcal M)$ of future-pointing curves such that $g(z)[\dot z,\dot z]=0, \ \text{a. e. on $[0,1]$}$. Consider the map $\Psi\colon \Omega_{p_0,q_0}(\mathcal M_0)\to \Lambda_{p,\gamma}(\mathcal M)$ defined by
\[\Psi(x)(s)=\big(x(s),t_0+\int_0^s F(x,\dot x)\de \nu\big).\]
Observe that $I$ is a smooth functional and $\Psi$ is differentiable at any regular curve $x$. Clearly we have that $J(z)=(I\circ\Psi)(x)=2E(x)$
and,
for any $\xi,\eta\in T_x\Omega_{p_0,q_0}(\mathcal M_0)$,
$\de\Psi(x)[\xi]$ is equal to
\beq\label{depsi}
\de \Psi(x)[\xi](s)=\left(\xi(s),\int_0^s(F_x(x,\dot x)[\xi]+F_y(x,\dot x)[\dot \xi])\de s\right),\eeq
hence $\de \Psi(x)$ is an injective map.
Since $x$ is a critical point of the length functional $x\mapsto\int_0^1F(x,\dot x)\de s$, we have that for any $\xi\in T_x\Omega_{p_0,q_0}(\mathcal M_0)$, $\de \Psi(x)[\xi](0)=\de \Psi(x)[\xi](1)=0$. Let now $U(s)=(U_0(s),\tau(s))$ be a vector field along  $z$ such that $U(0)=U(1)=0$ and $g(z)[\dot z,U]=0$. We  are going to show that $\de\Psi(x)[U_0]=U$ and hence $\de \Psi(x)$ is an isomorphism between the space of piecewise smooth vector  fields  along $x$ vanishing at the endpoints and the space of admissible variations for  $J$ (see Theorem~\ref{Uhl}).
Observe that $g(z)[U,\dot z]=0$ implies that
\[\tilde g_0(x)[U_0,\dot x]+\tilde g_0(x)[\delta(x),U_0]\dot t+ \tilde g_0(x)[\delta(x),\dot x]\tau-\tau\dot t=0.\]
By \eqref{fpnull} we get
\[\tilde g_0(x)[U_0,\dot x]+\tilde g_0(x)[\delta(x),U_0]F(x,\dot x)-\tau\sqrt{\alpha(x)[\dot x,\dot x]}=0,\]
where $\alpha(x)[\dot x,\dot x]=\tilde g_0(x)[\dot x,\dot x]+\tilde g_0(x)[\delta(x),\dot x]^2$.
Hence \[\tau=\frac{\tilde g_0(x)[U_0,\dot x]+\tilde g_0(x)[\delta(x),U_0]F(x,\dot x)}{\sqrt{\alpha(x)[\dot x,\dot x]}}.\]
From \eqref{depsi}, since $x$ is a geodesic for the metric $F$ and $U_0(0)=0$, the $t$ component of the vector field $\de\Psi(x)[U_0]$ is equal to $F_y(x,\dot x)[U_0(s)]$,  which is given by
\baln
\lefteqn{F_y(x,\dot x)[U_0(s)]=}&\\
&=\tilde g_0(x)[\delta(x),U_0]+\frac{\tilde g_0(x)[\delta(x),U_0]\tilde g_0(x)[\delta(x),\dot x]+\tilde g_0(x)[U_0,\dot x]}{\sqrt{\alpha(x)[\dot x,\dot x]}}=\tau(s).\ealn
Let  $\varphi=\varphi(r,s)\colon(-\eps,\eps)\times[0,1]\to \mathcal M$ be a variation defined by the admissible variational vector field $U=(U_0,\tau)$, and $\varphi_0=\varphi_0(r,s)\colon(-\eps,\eps)\times[0,1]\to \mathcal M_0$ be the one defined by $U_0$, we have that
\bml\label{ordine2}D^2 J(z)[U,U]=\frac{\de^2}{\de r^2}J(\varphi(r,\cdot))_{\big| r=0}\\
=\frac{\de^2}{\de r^2}I(\Psi(\varphi_0(r,\cdot)))_{\big| r=0}=
2\frac{\de^2}{\de r^2}E(\varphi_0(r,\cdot))_{\big| r=0}=2D^2 E(x)[U_0,U_0].\eml
From \eqref{ordine2}, by polarization, we get the equality between \eqref{hessJ} and the index form of the metric $F$ and then the equality \eqref{Morseindex}.
\end{proof}
\end{section}
\section*{Acknowledgments}
We wish to thank A. Abbondandolo for having drawn our attention to the paper \cite{LiLiLi05}.

\end{document}